\documentclass[11pt]{article}

\usepackage{amsmath,amsfonts,amssymb,a4, graphicx,bm}
\usepackage{amsfonts}
\usepackage{graphicx}
\usepackage{color}
\usepackage{color}
\usepackage{epsfig}
\usepackage{epstopdf}
\usepackage{tikz}
\usepackage{amstext}
\usetikzlibrary{arrows,decorations.pathmorphing,backgrounds,positioning,fit,petri}

\tikzset{nodde/.style={circle,draw=blue!50,fill=blue!20,inner sep=1.2pt}}

\begin{document}

\newtheorem{theorem}{Theorem}  [section]
\newtheorem{proposition}[theorem]{Proposition}
\newtheorem{lemma}[theorem]{Lemma}
\newtheorem{claim}[theorem]{Claim}
\newtheorem{cor}[theorem]{Corollary}
\newtheorem{fact}[theorem]{Fact}
\newtheorem{defn}[theorem]{Definition}
\newtheorem{conj}[theorem]{Conjecture}

\newcommand{\bproof}{\noindent{\bf Proof: }}
\newcommand{\eproof}{\hfill $\blacksquare$\\}
\newcommand{\eproofstate}{\hfill $\blacksquare$}
\newcommand{\de}{\mbox{def}}
\newcommand{\cee}{{\mathbb C}}
\newcommand{\real}{{\mathbb R}}
\newcommand{\R}{{\mathbb R}}
\newcommand{\rat}{{\mathbb Q}}
\newcommand{\zed}{{\mathbb Z}}
\newcommand{\calx}{{\cal X}}
\newcommand{\tran}{{{td}}}
\newcommand{\rank}{{\mbox{rank }}}
\newcommand{\sm}{\setminus}
\newcommand{\scrx}{{\mathcal X}}
\newcommand{\scrr}{{\mathcal R}}

\title{A necessary condition for generic rigidity of bar-and-joint frameworks in $d$-space}

\author{Hakan Guler\thanks{Department of Mathematics, Faculty of Arts \& Sciences, 
 Kastamonu University,
Kastamonu, Turkey. e-mail:
hakanguler19@gmail.com.}\\and\\Bill Jackson\thanks{School of Mathematical Sciences,
 Queen Mary University of London,
Mile End Road, London E1 4NS, England. e-mail:
b.jackson@qmul.ac.uk. }
}

\date{}

\maketitle

\begin{abstract}
A graph $G=(V,E)$ is {$d$-sparse} if each subset $X\subseteq V$ with
$|X|\geq d$ induces at most $d|X|-{{d+1}\choose{2}}$ edges in $G$.
 Maxwell showed in 1864 that a necessary condition for a generic bar-and-joint framework with at least $d+1$ vertices 
to be rigid in $\real^d$ is that $G$ should have a $d$-sparse subgraph with
$d|X|-{{d+1}\choose{2}}$ edges. This necessary condition is also sufficient when $d=1,2$ but not when $d\geq 3$. 
  Cheng and  Sitharam strengthened Maxwell's condition  by showing that {\em every} maximal $d$-sparse subgraph of 
  $G$ 
should have $d|X|-{{d+1}\choose{2}}$ edges when $d=3$. We extend their result to all $d\leq 11$.
\end{abstract}

\section{Introduction}
A {\em $d$-dimensional (bar-and-joint) framework} is a pair $(G,p)$
where $G=(V,E)$ is a graph and $p:V\to \real^d$. It is a long
standing open problem to determine when a given bar-and-joint
framework is {\em rigid} i.e.\ every continuous motion of the points
$p(v)$ which preserves the distances $\|p(u)-p(v)\|$ for all $uv\in
E$ must also preserve the distances $\|p(u)-p(v)\|$ for all $u,v\in
V$. It is not difficult to see that a 1-dimensional framework
$(G,p)$ is rigid if and only if the graph $G$ is connected. Abbot
\cite{A}  showed that the problem of determining rigidity
is NP-hard for all $d\geq 2$ but the  problem becomes more tractable if we assume that the
framework is {\em generic} i.e.\ there  are no algebraic dependencies
between the coordinates of the points $p(v)$, $v\in V$. 

Given a graph $G=(V,E)$, we can define a $|E|\times
d|V|$ matrix, the {\em $d$-dimensional rigidity matrix $R_d(G)$},
whose entries are linear combinations of indeterminates representing
the coordinates of the points $p(v)$,  in such a way that a generic framework $(G,p)$ with at least $d+1$ vertices  is
rigid if and only if the rank $r_d(G)$ of $R_d(G)$ is equal to
$d|V|-{{d+1}\choose{2}}$. This naturally gives rise to a matroid on
$E$, the {\em $d$-dimensional rigidity matroid ${\cal R}_d(G)$} in
which a set of edges $F\subseteq E$ is {\em independent}
if and only if the corresponding rows of $R_d(G)$ are linearly
independent. We refer the reader to \cite{W} for a precise
definition of the rigidity matrix, the rigidity matroid, and other
information on the topic of combinatorial rigidity.

Pollaczek-Geiringer \cite{PG} and subsequently Laman \cite{L} characterized when a 2-dimensional generic framework
is rigid (see also Lov\'asz and Yemini \cite{LY}). Their
characterization is based on the following concept. We say that a
graph $G=(V,E)$ is {\em $d$-sparse} if each  $X\subseteq V$ with $|X|\geq d+1$ induces at most $d|X|-{{d+1}\choose{2}}$
edges of $G$. Maxwell \cite{M} showed that being $d$-sparse is a
necessary condition for the rows of $R_d(G)$ 
to be linearly independent. Pollaczek-Geiringer  and Laman showed that that this
condition is also sufficient when $d=2$ and deduced that a
2-dimensional generic framework $(G,p)$ is rigid if and only if it
has a 2-sparse subgraph with $2|V|-3$ edges. Since every 
independent set of edges  in $\scrr_2(G)$ can be extended to a base of
 $\scrr_2(G)$, Laman's theorem implies that every
maximal 2-sparse subgraph of $G$ has the same number of edges.

It is known that the condition that $H$ is a $d$-sparse subgraph of $G$ is not
sufficient for  the edges of $H$ to
be  independent in $\scrr_d(G)$ when $d\geq 3$.  Indeed it is not even true
that all maximal $d$-sparse subgraphs of $G$ have the same number of
edges when $d\geq 3$. On the other hand, Cheng and Sitharam \cite{CS}
have shown that the number of edges in any maximal
$d$-sparse subgraph of $G$ does at least give an upper bound on
$r_d(G)$ when $d=3$. The purpose of this paper is to prove a result, Theorem
\ref{upperbound} below, which extends Cheng and
Sitharam's theorem to all values of $d\leq 11$.

\section{Sparse subgraphs}\label{sec:sparse}
Let $G=(V,E)$ be a graph and $d\geq 1$ be an integer. For
$X\subseteq V$ we use $E_G(X)$ to denote the set, and $i_G(X)$ the
number, of edges of $G$ joining pairs of vertices of $X$. We
simplify these to $E(X)$ and $i(X)$ when it is obvious to which
graph we are referring. We may rewrite the condition for $G$ to be
{$d$-sparse} as $i(X)\leq d|X|-{{d+1}\choose{2}}$ for all
$X\subseteq V$ with $|X|\geq d$. (Note that if $|X|\in \{d,d+1\}$
then we have $i(X)\leq {{|X|}\choose{2}}=d|X|-{{d+1}\choose{2}}$ and
the inequality holds trivially.) We will use the fact that the function $i:2^V\to \zed$ is {\em supermodular} i.e.\ $i(X)+i(Y)\leq i(X\cup Y)+i(X\cap Y)$ for all $X,Y\subseteq V$.

A subgraph $H=(U,F)$ of a $d$-sparse graph $G$ is {\em $d$-critical}
if either $|U|=2$ and $|F|=1$, or $|U|\geq d$ and
$|F|=d|X|-{{d+1}\choose{2}}$. The assumption that $G$ is $d$-sparse
implies that every $d$-critical subgraph of $G$ is an induced
subgraph. A {\em $d$-critical component} of $G$ is a $d$-critical
subgraph which is not properly contained in any other $d$-critical
subgraph of $G$.

\begin{lemma}\label{int}
Let $G=(V,E)$ be a $d$-sparse graph and
$H_1=(U_1,F_1),H_2=(U_2,F_2)$ be distinct $d$-critical components of
$G$. Then $|U_1\cap U_2|\leq d-1$ and, if equality holds, then
$i_G(U_1\cap U_2)={{d-1}\choose{2}}$.
\end{lemma}
\bproof Suppose that $|U_1\cap U_2|\geq d-1$. When $|U_1\cap
U_2|\geq d$ we have $i(U_1\cap U_2)\leq d|U_1\cap
U_2|-{{d+1}\choose{2}}$ since $G$ is $d$-sparse. When $|U_1\cap
U_2|=d-1$, we have $i(U_1\cap U_2)\leq {{d-1}\choose{2}}=d|U_1\cap
U_2|-{{d+1}\choose{2}}+1$ trivially. The maximality of $H_1,H_2$ and
the definition of a $d$-critical component imply that $|U_1|,|U_2|\geq
d$, and $d(|U_1|+|U_2|)-2{{d+1}\choose{2}}= i_G(U_1)+i_G(U_2)\leq
i_G(U_1\cup U_2)+i_G(U_1\cap U_2) \leq d|U_1\cup
U_2|-{{d+1}\choose{2}}-1+d|U_1\cap
U_2|-{{d+1}\choose{2}}+1=d(|U_1|+|U_2|)-2{{d+1}\choose{2}}.$
Equality must hold throughout. In particular we have  $i_G(U_1\cap
U_2)=d|U_1\cap U_2|-{{d+1}\choose{2}}+1$. This implies that $|U_1\cap U_2|= d-1$
and $i_G(U_1\cap U_2)={{d-1}\choose{2}}$. \eproof

Let $k,t$ be non-negative integers, $G=(V,E)$ be a graph and $\cal X$
be a family of subsets of $V$.
We say that $\scrx$ is 
 {\em $t$-thin} if every pair of sets in $\calx$
intersect in at most $t$ vertices.
A {\em $k$-hinge} of $\cal X$ is a
set of $k$ vertices which lie in the intersection of at least two
sets in $\cal X$. A $k$-hinge $U$ of $\calx$ is {\em closed in $G$}
if $G[U]$ is a complete graph. We use $\Theta_k(\calx)$
to denote the set of all
$k$-hinges
of $\cal X$. For $U\in
\Theta_k(\calx)$, let $d_\calx(U)$ denote the number of sets in
$\cal X$ which contain $U$. Note that if $G$ is $t$-thin then
$\Theta_k(\calx)=\emptyset$ for all $k\geq t+1$. Note also that
$\Theta_0(\calx)=\{\emptyset\}$ and $d_\calx(\emptyset)=|\calx|$.

\begin{lemma}\label{prefixedhinge}
Let $H=(V,E)$ be a $d$-sparse graph, $\cal X$ be a family of subsets of $V$ such that $H[V_i]$ is $d$-critical for all $V_i\in \scrx$, and $W\in \Theta_k(\calx)$ for some $0\leq k\leq d-1$. Suppose
that
$|V_i|\geq d$ for all $V_i\in \scrx$ with $W\subseteq V_i$. Then
$$
(d-k)\sum_{\substack{U\in\Theta_{k+1}(\calx)\\[0.7mm]W\subset U}}(d_\calx(U)-1)-
\sum_{\substack{U\in\Theta_{k+2}(\calx)\\[0.7mm]W\subset U}}(d_\calx(U)-1)\leq{{d+1-k}\choose{2}}(d_\calx(W)-1)
\,.$$
\end{lemma}
\bproof Let $d_\calx(W)=t$ and let $V_1,V_2,\ldots,V_t$ be the
sets in $\scrx$ which contain $W$. Let $H_i=(V_i,E_i)=H[V_i]$
for $1\leq i\leq t$. Let $H'=\bigcup_{i=1}^tH_i$ and put
$H'=(V',E')$. Then
\begin{equation}\label{fixedhingeeq1}
 |V'|=\sum_{i=1}^t|V_i|-k(t-1)-\sum_{\substack{U\in\Theta_{k+1}(\calx)\\[0.7mm]W\subset
 U}}(d_\calx(U)-1)
\end{equation}
since, for $v\in V'$, if $v\in W$ then $v$ is counted $t$ times in
$\sum_{i=1}^t|V_i|$, if $v\in U\setminus W$ for some
$U\in\Theta_{k+1}$ with $W\subset U$ then $v$ is counted
$d_\calx(U)$ times in $\sum_{i=1}^t|V_i|$, and all other vertices of
$V'$ are counted exactly once in $\sum_{i=1}^t|V_i|$.

Similarly,
\begin{equation}\label{fixedhingeeq2}
 |E'|\geq\sum_{i=1}^t|E_i|-{{k}\choose{2}}(t-1)-k\sum_{\substack{U\in\Theta_{k+1}(\calx)\\[0.7mm]W\subset
 U}}(d_\calx(U)-1)-\sum_{\substack{U\in\Theta_{k+2}(\calx)\\[0.7mm]W\subset
 U}}(d_\calx(U)-1)
\end{equation}
since, for $e=xy\in E'$: if $x,y\in W$ then $e$ is counted $t$ times
in $\sum_{i=1}^t|E_i|$ and there are at most ${{k}\choose{2}}$ such
edges; if $x\in W$ and $y\in U\setminus W$ for some
$U\in\Theta_{k+1}$ with $W\subset U$ then $e$ is counted
$d_\calx(U)$ times in $\sum_{i=1}^t|E_i|$ and for each such $y$
there are at most $k$ choices for $x$; if $x,y\in U\setminus W$ for
some $U\in\Theta_{k+2}$ with $W\subset U$ then $e$ is counted
$d_\calx(U)$ times in $\sum_{i=1}^t|E_i|$, and all other edges of
$E'$ are counted exactly once in $\sum_{i=1}^t|E_i|$.

Since $H'\subseteq H$, $H'$ is $d$-sparse
Hence $|E'|\leq
d|V'|-{{d+1}\choose{2}}$. We may substitute equations
(\ref{fixedhingeeq1}) and (\ref{fixedhingeeq2}) into this inequality
and use the fact that $|E_i|=d|V_i|-{{d+1}\choose{2}}$ for all
$1\leq i\leq t$ to obtain
\begin{eqnarray*}
(d-k)\sum_{\substack{U\in\Theta_{k+1}(\calx)\\[0.7mm]W\subset
U}}(d_\calx(U)-1)&-&
\sum_{\substack{U\in\Theta_{k+2}(\calx)\\[0.7mm]W\subset
U}}(d_\calx(U)-1)\\
&\leq&\left[{{d+1}\choose{2}}+{{k}\choose{2}}-dk\right](t-1)\\
&=& {{d+1-k}\choose{2}}(t-1).
\end{eqnarray*}
\eproof

\begin{lemma}\label{fixedhinge}
Let $H=(V,E)$ be a $d$-sparse graph, $\cal X$ be a
family of subsets of $V$ such that $H[V_i]$ is $d$-critical
and
$|V_i|\geq d$ for all $V_i\in \scrx$.
 Put
$a_k=\sum_{U\in\Theta_k(\calx)}(d_\calx(U)-1)$ for  $0\leq k\leq d$.
Then
for all $0\leq k\leq d-2$ we have:\\[2mm]
(a) $(d-k)(k+1)a_{k+1}-{{k+2}\choose{2}}a_{k+2}
\leq{{d+1-k}\choose{2}}a_{k}$;
\\[2mm]
(b) $(d-k)a_{k+1}- (k+1)a_{k+2} \leq{{d+1}\choose{k+2}}(|\calx|-1)$;
\\[2mm]
(c) if $\scrx$ is $(d-1)$-thin, $d(d-k)a_{k+1} \leq(k+2)(d-k-1){{d+1}\choose{k+2}}(|\calx|-1)$.
\end{lemma}
\bproof Part (a) follows by summing the inequality in Lemma
\ref{prefixedhinge} over all $W\in \Theta_k$, and using the facts
that
$$\sum_{W\in \Theta_k(\calx)} \;\sum_{\substack{U\in\Theta_{k+1}(\calx)\\[0.7mm]W\subset
U}}(d_\calx(U)-1)=(k+1)\sum_{U\in\Theta_{k+1}(\calx)}(d_\calx(U)-1)=(k+1)a_{k+1}$$
and
$$\sum_{W\in \Theta_k(\calx)} \;\sum_{\substack{U\in\Theta_{k+2}(\calx)\\[0.7mm]W\subset
U}}(d_\calx(U)-1)={{k+2}\choose{2}}\sum_{U\in\Theta_{k+2}(\calx)}(d_\calx(U)-1)={{k+2}\choose{2}}a_{k+2}\,.$$

We prove (b) by induction on $k$. When $k=0$, (b) follows by putting
$k=0$ in (a), and using the fact that $a_0=|\scrx |-1$. Hence suppose
that $k\geq 1$. Then (a) gives
\begin{equation}\label{fixedhingeeq3}
2(d-k)a_{k+1}-2(k+1)a_{k+2}\leq\frac{(d-k+1)(d-k)}{k+1}\,a_k-ka_{k+2}\,.
\end{equation}
We may also use (a) to obtain
\begin{equation}\label{fixedhingeeq4}
ka_{k+2}\geq\frac{k(d-k)}{k+2}\left(2a_{k+1}-\frac{d-k+1}{k+1}a_{k}\right)\,.
\end{equation}
Substituting (\ref{fixedhingeeq4}) into (\ref{fixedhingeeq3}) and
using induction we obtain
\begin{eqnarray*}
(d-k)a_{k+1}-(k+1)a_{k+2}&\leq&\mbox{$\frac{d-k}{k+2}\,[(d-k+1)a_k-ka_{k+1}]$}\\
&\leq&\mbox{$\frac{d-k}{k+2}\,{{d+1}\choose{k+1}}\,(|\calx|-1)$}\\
&=&\mbox{${{d+1}\choose{k+2}}\,(|\calx|-1)$}\,.
\end{eqnarray*}

We prove (c) by induction on $d-k$. When $d-k=2$, (c) follows by
putting $k=d-2$ in (b)  and using the fact that $a_d=0$
since $\scrx$ is $(d-1)$-thin.
Hence suppose that $d-k\geq 3$. Then (b) gives
$$
d(d-k)a_{k+1}\leq\mbox{$d{{d+1}\choose{k+2}}\,(|\calx|-1)+d(k+1)a_{k+2}$}\,.
$$
We may now apply induction to $a_{k+2}$ to obtain
\begin{eqnarray*}
d(d-k)a_{k+1}&\leq&\mbox{$[d{{d+1}\choose{k+2}}+\frac{(k+1)(k+3)(d-k-2)}{d-k-1}{{d+1}\choose{k+3}}]
\,(|\calx|-1)$}\\[1mm]
&=&\mbox{$(k+2)(d-k-1){{d+1}\choose{k+2}} \,(|\calx|-1)$}\,.
\end{eqnarray*}
\eproof

\begin{theorem}\label{boundedhinges}
Let $H=(V,E)$ be a $d$-sparse graph, $\cal X$ be a $(d-1)$-thin
family of subsets of $V$ such that $H[X]$ is $d$-critical and $|X|\geq d$
 for all $X\in \scrx$.
For each  $X\in \scrx$  let $\theta_k(X)$
be the number of $k$-hinges of $\calx$ contained in $X$.
Then:\\[1mm]
(a) $\theta_{1}(X)\leq 2d-1$ for some $X\in \scrx$;
\\[1mm]
(b) $\theta_{2}(X)\leq (d-2)(d+1)-1$ for some $X\in \scrx$;
\\[1mm]
(c) $\theta_{d-1}(X)\leq d$ for some  $X\in \scrx$.
\end{theorem}
\bproof 

We first prove (a). Putting $k=0$ in Lemma \ref{fixedhinge}(c) we
obtain
\begin{equation}\label{boundedhingeseq1}
d\sum_{U\in\Theta_{1}(\calx)}(d_\calx(U)-1)
\leq(d-1)(d+1)(|\calx|-1)\,.
\end{equation}
Since $d_\calx(U)\geq 2$ for all $U\in\Theta_{1}(\calx)$ we have
$d_\calx(U)-1\geq d_\calx(U)/2$ and hence (\ref{boundedhingeseq1})
gives
$$\sum_{U\in\Theta_{1}(\calx)}d_\calx(U)
<2d\,|\calx|\,.
$$
This tells us that the average number of $1$-hinges in a set in $\scrx$ is strictly less than $2d$.

We next prove (b). Putting $k=1$ in Lemma \ref{fixedhinge}(c) we
obtain
\begin{equation}\label{boundedhingeseq1.5}
\sum_{U\in\Theta_{2}(\calx)}(d_\calx(U)-1)
\leq(d-2)(d+1)(|\calx|-1)/2\,.
\end{equation}
We can now proceed as in (a).

Finally we prove (c). Putting $k=d-2$ in Lemma \ref{fixedhinge}(c)
gives
\begin{equation}\label{boundedhingeseq3}
2\sum_{U\in\Theta_{d-1}(\calx)}(d_\calx(U)-1) \leq(d+1)(|\calx|-1)\,.
\end{equation}
We can now proceed as in (a).
\eproof

The bounds given in Theorem \ref{boundedhinges} (a), (b) are close to being best possible. To see this consider the graph $H=H_1\cup H_2\cup\ldots\cup H_m$ where $H_i=(V_i,E_i)$ is $d$-critical, $H_i\cap H_{j}=K_{d-1}$ for $i-j\equiv \pm 1 \mod m$ and otherwise $H_i\cap H_{j}=\emptyset$. Then $H$ is $d$-sparse when $m$ is sufficiently large,  $\scrx=\{V_1,V_2,\ldots, V_m\}$ is $(d-1)$-thin and we have $\theta_1(V_i)=2d-2$ and $\theta_2(V_i)=(d-1)(d-2)$ for all $V_i\in \scrx$. We do not know whether (c) is close to best possible for large $d$. It is conceivable that there  always exists a set $X\in \scrx$ with $\theta_{d-1}(X)\leq 2$.

\section{Main result}

In order to prove our main theorem we will need the following result from \cite{GJN}.

\begin{lemma}\label{lem:circuit}
Let $G=(V,E)$ be a graph such that $E$ is a non-rigid circuit in $\scrr_d(G)$. Then 
$|E|\geq d(d+9)/2$. \eproofstate
\end{lemma}

Let  $G=(V,E)$ be a graph and $\cal X$
be a family of subsets of $V$.
We say that $\cal X$ is a {\em cover}
of $G$ if every set in $\calx$ contains at least two vertices, and
every edge of $G$ is induced by at least one set in $\cal X$. 

\begin{lemma}\label{cover}
Let $G=(V,E)$ be a graph, $H=(V,F)$ be a maximal $d$-sparse subgraph
of $G$, and  $H_1,H_2,\ldots,H_m$ be the  $d$-critical components of
$H$. Let $X_i$ be the vertex set of $H_i$ for  $1\leq i \leq m$.
Then ${\cal X}=\{X_1, X_2,\ldots, X_m\}$ is a $(d-1)$-thin cover of
$G$ and each $(d-1)$-hinge of $\calx$ is closed in $H$.
\end{lemma}
\bproof The definition of a $d$-critical subgraph implies that each
$H_i$ has at least two vertices and that every edge of $H$ belongs
to at least one $H_i$. Thus $\calx$ is a cover of $H$. To see that
$\calx$ also covers $G$ we choose $e=uv\in E\sm F$. The maximality
of $H$ implies that $H+e$ is not $d$-sparse. Hence $\{u,v\}$ is
contained in some $d$-critical subgraph of $H$. Thus $\calx$ also
covers $G$. The facts that $\cal X$ is $(d-1)$-thin and that each
$(d-1)$-hinge of $\calx$ is closed follow from Lemma \ref{int}.
\eproof

We refer to the $(d-1)$-thin cover of $G$ described in Lemma
\ref{cover} as the {\em $H$-critical cover} of $G$. 
Note that the definition of a
$d$-critical set implies that each set in the $H$-critical cover has
size two or has size at least $d$.

\begin{theorem}\label{upperbound}
Let $G=(V,E)$ be a graph, $d\leq 11$ be an integer and $H=(V,F)$ be a
maximal $d$-sparse subgraph of $G$.  Then $r_d(G)\leq |F|$.
\end{theorem}
\bproof We proceed by contradiction. Suppose the theorem is false
and choose a counterexample $(G,H)$ such that $|E|$ is as small as
possible. Let $H_1,H_2,\ldots,H_m$ be the $d$-critical components of
$H$ where $H_i=(V_i,F_i)$ for $1\leq i\leq m$. Then
$\calx_{0}=\{V_1,V_2,\ldots,V_m\}$ is the $H$-critical cover of $G$. 

Choose a cover $\scrx$ of $G$ such that $\scrx\subseteq \scrx_0$ and $|\scrx|$ is as small as 
possible. Note that $\scrx_0$, and hence also $\scrx$, are $(d-1)$-thin. For each $V_i\in \scrx$, let 
$F_i^*$ be the set of all edges $uv\in F_i$ such that $\{u,v\}$ is a 2-hinge of $\scrx$,  and let $E_i$ be the
set of edges of $G$ induced by $V_i$.

\begin{claim}\label{c1}
If $e=uv\in E$ satisfies $r_d(G)=r_d(G-e)$, then $\{ u,v\}$ is a 2-hinge of $\scrx$.
\end{claim}
\bproof First suppose that $e\in E\setminus F$. Since $H$ is a maximal $d$-sparse
subgraph of $G-e$, the minimality of $|E|$ gives $r_d(G-e)\leq |F|$.  Since $r_d(G)=r_d(G-e)$ this gives
a contradiction.

Thus we can assume that $e\in F$. Let $d_{\scrx}(e)$ be the number of $V_i\in \scrx$ such
that $e\in F_i$. Since $H-e$ is a $d$-sparse subgraph
of $G-e$, we may choose a maximal $d$-sparse subgraph $H'=(V,F')$ of $G-e$ which contains $H-e$.
Let $V_i\in \scrx$. 
If $e\notin F_i$, then no edge of $E_i\setminus F_i$ can be in $F'$, since
$H_i$ is $d$-critical.  On the other hand,  if $e\in F_i$, then at most one edge of
$E_i\setminus F_i$ can be in $F'$, since $|F_i-e|=d|V_i|-{{d+1}\choose 2 }-1$.
These observations imply that $|F'|\leq |F|-1+d_{\scrx}(e)$. By the minimality of $|E|$
we have $r_d(G-e)\leq |F'|$, and hence $r_d(G)=r_d(G-e)\leq |F|-1+d_{\scrx}(e)$. Combining this with
$r_d(G)>|F|$ gives $d_{\scrx}(e)\geq 2$.\eproof

We next show  that $F_i^*$ is dependent in $\scrr_d(G)$ for all $V_i\in \scrx$. Suppose this is
not the case. Then $E_i$ is independent in $\scrr_d(G)$ by Claim \ref{c1}.
Thus $E_i$ can have at most $d|V_i|-{{d+1}\choose 2 }$ edges. Since  $H_i$
is $d$-critical, this gives $E_i=F_i$. 
The  minimality of $\scrx$ implies that $F_i\neq F_i^*$ and hence we may choose  an edge $e\in F_i\setminus F_i^*$. 
Since $F_i=E_i$, all edges
of $G-e$ which are induced by $V_i$ are in $H-e$. Since each $V_j\in \scrx-V_i$ induce a  $d$-critical subgraph of $H-e$, we conclude that $H-e$ is a maximal $d$-sparse subgraph of $G-e$. The minimality of $|E|$ now gives $r_d(G-e)\leq |F-e|=|F|-1$.
Since $e\not\in F_i^*$,  Claim \ref{c1} gives $r_d(G-e)=r_d(G)-1$. Hence $r_d(G)=r_d(G-e)+1\leq |F|$.
This contradicts the choice of $G$ and implies that $F_i^*$ is dependent in $\scrr_d(G)$ for all $V_i\in \scrx$.

By  Theorem \ref{boundedhinges}(b) we may choose $V_i\in \scrx$ such that $|F_i^*|\leq (d-2)(d+1)-1$. Since $F_i^*$ is dependent in $\scrr_d(G)$, it contains a
circuit of ${\cal R}_d(G)$. This circuit cannot be rigid, since $H$
is $d$-sparse. Lemma \ref{lem:circuit} now gives  $\frac{d^2+9d}{2} \leq |F_i^*|\leq (d-2)(d+1)-1
$ which implies that $d\geq 12$.
\eproof

We have the following immediate corollary.

\begin{cor}\label{cor:rigid}
Let $d\leq 11$ be an integer and $G=(V,E)$ be a graph with $|V|\geq d+1$. If $G$ is generically rigid in $\R^d$ then every maximal $d$-sparse subgraph of $G$ has $d|V|-{{d+1}\choose{2}}$ edges. \eproofstate
\end{cor}

\section{Closing remarks}
1. Given a graph $G$, let $s_d(G)$ be the minimum number of edges in a
maximal $d$-sparse subgraph of $G$. Theorem \ref{upperbound} tells
us that $r_d(G)\leq s_d(G)$ when $d\leq 11$. 
We can use the
following operation to
construct graphs for which strict inequality holds. Given two graphs $G_1=(V_1,E_1)$ and
$G_2=(V_2,E_2)$ with $V_1\cap V_2=\{u,v\}$ and $E_1\cap E_2=\{uv\}$,
we refer to the graph $G=G_1\cup G_2$ as the {\em parallel
connection of $G_1$ and $G_2$ along the edge $uv$}.

The graph $G$ obtained by taking the parallel connection of two
copies of $K_5$ along an edge $uv$ and then deleting $uv$, is
3-sparse and is not rigid in $\real^3$. Hence
$s_3(G)=|E(G)|=18>17=r_3(G)$. On the other hand we may improve the
upper bound on $r_3(G)$ in this example by considering the graph
$H=G+uv$. A maximal 3-sparse subgraph of $H$ which contains $uv$
has 17 edges. Thus we have $17=r_3(G)\leq r_3(H)\leq s_3(H)=17$.

More generally, for any graph $G$, let $s_d^*(G)=\min\{s_d(H)\;:\;G\subseteq H\}$.
Then $r_d(G)\leq s_d^*(G)$ for all $d\leq 11$. 
 The following example shows that strict
inequality can also hold in this inequality. 
Let $G$ be obtained
from $K_5$ by taking parallel connections with 10 different $K_5$'s
along each of the edges of the original $K_5$. We have $r_3(G)=89$.
On the other hand,
 $s_3(G)=90$ (obtained by taking a maximal 3-sparse
subgraph which contains nine of the edges of the original $K_5$).
Furthermore we have $s_3(H)\geq r_3(H)> r_3(G)$ for all graphs
$H$ which properly contain $G$. Thus  $s_3^*(G)=90>r_3(G)$.

\medskip\noindent
2. For fixed $d$, we can use network flow algorithms to test whether a
graph is $d$-sparse in polynomial time, see for example \cite{BJ}.
This means we can greedily construct a maximal $d$-sparse subgraph
$H$ of a graph $G$ in polynomial time and hence obtain an upper
bound on $r_d(G)$ when $d\leq 11$ via Theorem \ref{upperbound}. 
We do not know
whether $s_d(G)$ or $s_d^*(G)$ can be determined in polynomial time.

\medskip\noindent
3. We believe that the conclusion of Theorem \ref{upperbound} should be valid for all $d$. However the graph $G$ given in the example at the end of Section \ref{sec:sparse} shows that
our proof technique will not give this: $G$ is $d$-sparse and we have $\theta_2(V_i)=(d-1)(d-2)$ for all $V_i$ in the $G$-critical cover of $G$. 
On the other hand, the lower bound on the number of edges in a non-rigid circuit in $\scrr_d(G)$ given by Lemma \ref{lem:circuit} is $\frac{d(d+9)}{2}$, so we cannot use it to deduce that
the set of 2-hinges in some $G[V_i]$ is $\scrr_d$-independent when $d\geq 15$. One way to get round this problem would be to show that the $d$-critical components in a $d$-sparse graph
form a cover which is `iteratively independent' i.e.\ we can order the vertex sets of these components as $V_1,V_2,\ldots V_m$ such that the set of 2-hinges of
$\{V_1,V_2,\ldots, V_i\}$ which belong to $V_i$ is $\scrr_d$-independent for all $2\leq i\leq m$. We refer the reader to \cite{JJiterative} for more information on iteratively
independent covers.

\subsubsection*{Acknowledgement} The authors would like to thank Meera Sitharam
for helpful conversations on this topic and the London Mathematical Society for providing partial financial support through a scheme 5 grant.
The first author would  also like to thank the Ministry of National Education of Turkey for PhD funding through a YLSY grant.

\end{document}